\documentclass[leqno,12pt]{article}
\usepackage{amsmath,amssymb,amscd,amsthm,latexsym,amsfonts, epsfig, graphicx}

\begin{document}

\numberwithin{equation}{section} 

\newtheorem{theorem}{Theorem}[section]
\newtheorem{question}{Question}
\newtheorem{conjecture}[theorem]{Conjecture} 
\newtheorem{lemma}[theorem]{Lemma}
\newtheorem*{claim}{Claim}
\newtheorem{corollary}[theorem]{Corollary}
\newtheorem{proposition}[theorem]{Proposition}

\theoremstyle{remark}
\newtheorem*{remark}{Remark}
\newtheorem*{remarks}{Remarks}
\newtheorem*{merci}{Acknowledgements}
\newtheorem*{defi}{Definition}

\newcommand{\dv}{\operatorname{div}}
\newcommand{\R}{\operatorname{Re}}
\newcommand{\supp}{\operatorname{supp}}
\newcommand{\dist}{\operatorname{dist}}
\newcommand{\Lip}{\operatorname{Lip}}
\newcommand{\diam}{\operatorname{diam}}
\newcommand{\epi}{\operatorname{Epi}}

\newcommand{\NN}{\mathbb{N}}
\newcommand{\RR}{\mathbb{R}}
\newcommand{\ZZ}{\mathbb{Z}}
\newcommand{\QQ}{\mathbb{Q}}
\newcommand{\CC}{\mathbb{C}}
\newcommand{\KK}{\mathbb{K}}
\newcommand{\rn}{\RR^n}
\newcommand{\srn}{{\scriptscriptstyle \RR}^n}
\newcommand{\rplus}{\RR_+}
\newcommand{\rplusbar}{\ol{\RR_+}}

\newcommand{\ep}{\varepsilon}
\newcommand{\si}{\sigma}
\newcommand{\dxdtovert}{\frac{dxdt}{t}}
\newcommand{\dtovert}{\frac{dt}{t}}
\newcommand{\comp}{{}^{\textnormal{c}}}
\newcommand{\lims}{\ol\lim}
\newcommand{\limi}{\ul\lim}
\newcommand{\adelta}{\calA(L)}
\newcommand{\om}{\Omega}
\newcommand{\pom}{\partial\om}
\newcommand{\dive}{\mathrm{div}\,}
\newcommand{\divagrad}{\,\dive\, (A \nabla)}
\newcommand{\kloc}{(K$_{\mathrm{loc}}$)}
\newcommand{\gloc}{(G$_{\mathrm{loc}}$)}
\newcommand{\dloc}{(D$_{\mathrm{loc}}$)}
\newcommand{\CZ}{Calder\'on-Zygmund }

\newcommand{\ul}{\underline}
\newcommand{\ol}{\overline}
\newcommand{\exist}{\;\exists\,}
\newcommand{\foral}{\;\forall\,}

\newcommand{\calA}{\mathcal{A}}
\newcommand{\calB}{\mathcal{B}}
\newcommand{\calC}{\mathcal{C}}
\newcommand{\calD}{\mathcal{D}}
\newcommand{\calE}{\mathcal{E}}
\newcommand{\calF}{\mathcal{F}}
\newcommand{\calG}{\mathcal{G}} 
\newcommand{\calH}{\mathcal{H}} 
\newcommand{\calI}{\mathcal{I}} 
\newcommand{\calJ}{\mathcal{J}} 
\newcommand{\calK}{\mathcal{K}} 
\newcommand{\calL}{\mathcal{L}}
\newcommand{\calM}{\mathcal{M}}
\newcommand{\calN}{\mathcal{N}}
\newcommand{\calO}{\mathcal{O}}
\newcommand{\calP}{\mathcal{P}}
\newcommand{\calQ}{\mathcal{Q}}
\newcommand{\calR}{\mathcal{R}}
\newcommand{\calS}{\mathcal{S}}
\newcommand{\calT}{\mathcal{T}}
\newcommand{\calU}{\mathcal{U}}
\newcommand{\calV}{\mathcal{V}}
\newcommand{\calW}{\mathcal{W}}
\newcommand{\calX}{\mathcal{X}}
\newcommand{\calY}{\mathcal{Y}}
\newcommand{\calZ}{\mathcal{Z}}

\title{
On the Calder\'on-Zygmund lemma for Sobolev functions}

\author{Pascal Auscher~\footnote{Universit\'e  Paris-Sud and CNRS UMR 
8628, Department de Math\'ematiques, 91405 Orsay Cedex (France) Email: 
\texttt{pascal.auscher@math.u-psud.fr}}}  

\date{october 16, 2008}
   
\maketitle

\begin{abstract} We correct an inaccuracy in the proof of a result in \cite{Aus1}.
\end{abstract}

2000 MSC: 42B20, 46E35

Key words: Calder\'on-Zygmund
decomposition; Sobolev spaces.

We recall the lemma.

\begin{lemma}\label{lemmaCZD} Let $n\ge 1$,  $1\le p\le \infty$ and $f\in \calD'(\RR^n)$ be such that 
$\|\nabla f\|_p <\infty$. Let $\alpha>0$. Then, one can find a collection of cubes $(Q_i)$, functions $g$ and $b_i$  such that 
\begin{equation}\label{eqcsds1}
f= g+\sum_i b_i  \end{equation}
and the following properties hold:
\begin{equation}\label{eqcsds2}
\|\nabla g\|_\infty \le C\alpha, \end{equation}
\begin{equation}\label{eqcsds3}
b_i \in W_0^{1,p}(Q_i)\ \text{and} \ \int_{Q_i} |\nabla b_i|^p \le C\alpha^p |Q_i|, \end{equation}
 \begin{equation}\label{eqcsds4}
\sum_i |Q_i| \le C\alpha^{-p} \int_{\RR^n} |\nabla f|^p , \end{equation}
\begin{equation}\label{eqcsds5}
\sum_i {\bf 1}_{Q_i} \le N, \end{equation}
where $C$ and
 $N$ depend only on dimension and $p$.
\end{lemma}

This lemma was first  stated in \cite{Aus1} in $\RR^n$. Then it appears in various forms and extensions in 
\cite{Aus-mem}  (same proof), \cite{AC} (same proof on manifolds), \cite{AM} (on $\RR^n$ but with a doubling weight), B. Ben Ali's PhD thesis \cite{Be} and \cite{AB} (The Sobolev space is modified to adapt to Schr\"odinger operators), N. Badr's  PhD thesis \cite{Ba} and \cite{Ba1, Ba2}  (used toward interpolation of Sobolev spaces on manifolds and measured metric spaces) and in \cite{BR} (Sobolev spaces on graphs). The same inaccuracy can be corrected everywhere as below. 
The proof of the generalisation to higher order Sobolev spaces in \cite{Aus1}  can also be corrected  with similar ideas.

The second equation tells that $g$ is in fact Lipschitz continuous. There is a direct proof of this fact in N. Badr's thesis \cite{Ba}. The proof proposed in \cite{Aus1} is as follows:

Define $b_i=(f-c_i)\calX_i$ where $c_i$ are appropriate numbers and $(\calX_i)$ forms a smooth partition of unity of $\Omega= \cup Q_i$ subordinated to the cubes $(\frac 1 2 Q_i)$ with support of $\calX_i$ contained in $Q_i$. For example,  the choice $c_i=f(x_i)$ for some well chosen $x_i$ or $c_i= m_{Q_i} f$, the mean of $f$ over the cube $Q_i$, ensures that $\sum_i |b_i| \ell_i^{-1}$ is locally integrable  ($\ell_i$ is the length of $Q_i$) and that $\sum_i b_i$ is a distribution on $\RR^n$. Then   $g$ defined  as $g= f- \sum_i b_i$ is a distribution on $\RR^n$. 
Its gradient  $\nabla g$ can be calculated as  $\nabla g= (\nabla f ){\bf 1}_F + h$ in the sense of distributions (on $\RR^n$) with 
$h= \sum_i c_i \nabla\calX_i$. It is then a consequence of the construction of the set $F=\Omega^c$ that $|\nabla f|$ is bounded on $F$  by $\alpha$ and then it is  shown that $|h|$ is bounded by $C\alpha$ , which implies the boundedness of $|\nabla g|$.

Everything is correct in the argument above BUT the representation of $h$. The series $\sum_i c_i \nabla\calX_i$, viewed as the distributional derivative on $\RR^n$ of $\sum_i c_i \calX_i$, may not be a measurable function (section) on $\RR^n$. For example, if  
$c_i=1$ for all $i$, then $\sum_i \nabla\calX_i = \nabla {\bf 1}_\Omega$  is  a non-zero distribution supported on the boundary of $\Omega$ (a measure if $\Omega$ has locally finite perimeter). One needs to renormalize the series to make it converge in the distribution sense.

Here we give  correct renormalizations of   $h$. A first one is obtained right away from differentiation of $g$:
$$
h= - \sum_i (f-c_i)\nabla \calX_i .
$$
The convergence in the distributional sense in $\RR^n$ is in fact hidden of \cite{Aus1}.

A second one is
$$
h= - \sum_m\left(\sum_i (c_m-c_i)\nabla \calX_i \right)\calX_m.
$$
This representation converges in the distributional sense in $\RR^n$ and can be shown to be a bounded function.

Let us show how to obtain the second representation in the sense of distributions. Then the proof of boundedness is as in \cite{Aus1}. Take a test function $\phi$ in $\RR^n$. Then by definition the distribution $\sum_i \nabla b_i$ tested against $\phi$ is given by 
$$
\sum_i  \int \nabla b_i \ \phi.
$$
To compute this, we take a finite subset $J$ of the set $I$ of indices $i$ and we have to pass to the limit in the sum restricted to $J$ as $J\uparrow I$.  Because now the sum is finite, and all functions have support in the set $\Omega$, we can introduce $\sum_m \calX_m= {\bf 1_\Omega}$. We have
$$
\sum_{i\in J}  \int \nabla b_i \ \phi = \sum_m \sum_{i\in J}  \int \nabla b_i \  \calX_m \ \phi.
$$
Now recall that $b_i= (f-c_i)\calX_i$ . Call $I_m$ the set of indices such that the support of $\calX_i$ meets the support of $\calX_m$. By property of the Whitney cubes, $I_m$ is a finite set with bounded cardinal. Hence we can write
$$
\sum_m \sum_{i\in J}  \int \nabla b_i \  \calX_m \ \phi = \sum_m \sum_{i\in J\cap I_m}  \int  \nabla f \ \calX_i \  \calX_m \ \phi + \sum_m \sum_{i\in J\cap I_m}  \int (f-c_i)\nabla \calX_i \  \calX_m \ \phi.
$$
It is clear that the first sum in the RHS converges to $\int_\Omega \nabla f \ \phi$ as $J\uparrow I$. 
As for the second it is equal to 
$$
\sum_m \sum_{i\in J\cap I_m}  \int  (c_m-c_i)\nabla \calX_i \  \calX_m \ \phi
+ \sum_m \sum_{i\in J\cap I_m}  \int (f-c_m)\nabla \calX_i \  \calX_m \ \phi
$$
As one can show (with the argument in \cite{Aus1}) that 
$$
\sum_m \sum_{i}   |c_m-c_i||\nabla \calX_i| | \calX_m| \le C\alpha
$$
one has  
$$
\lim_{J\uparrow I} \sum_m \sum_{i\in J\cap I_m}  \int  (c_m-c_i)\nabla \calX_i \  \calX_m \ \phi =- \int h \ \phi
$$
where $h$ is defined above. 

Finally, write  
$$\sum_m \sum_{i\in J\cap I_m}  \int  (f-c_m)\nabla \calX_i \  \calX_m \ \phi
=
\sum_m \int  b_m R_{m, J} \ \phi
$$
with 
$$
R_{m, J}= \sum_{i\in J\cap I_m}   \nabla \calX_i .
$$
By construction of the $\calX_i$ and properties of Whitney cubes, 
$$
\sum_{i\in I_m}   |\nabla \calX_i | \le C\ell_m^{-1}
$$
where $\ell_m$ is the length of $Q_m$,
and on the support of $\calX_m$
$$
\sum_{i\in I_m}  \nabla  \calX_i = \sum_{i\in I}   \nabla \calX_i =0.
$$
As  $\sum_i |b_m| \ell_m^{-1}$ has been shown to be locally integrable, one can conclude by the Lebesgue dominated convergence theorem that 
$$\lim_{J\uparrow I} \sum_m \int  b_m R_{m, J} \ \phi = 0.
$$

\end{document}